# A Note on the Primary Decomposition of $k-$Ideals in Semirings

Ram Parkash Sharma[1,a], Ricah Sharma[1,b], S. Kar[2,c] and Madhu[1,*]

[1]Department of Mathematics and Statistics, Himachal Pradesh University, Summer Hill, Shimla-171005, India

[2]Department of Mathematics, Jadavpur University, West Bengal, Kolkata-700032, India

*Corresponding author E-mail address: mpatial.math@gmail.com

a. rp_math_hpu@yahoo.com

b. richasharma567@yahoo.com

c. karsukhendu@yahoo.co.in

ABSTRACT: We establish the primary decomposition and uniqueness of primary decomposition for $k-$ideals in commutative Noetherian semirings.

## 1. Introduction

Ideals play an important role in both semiring theory and ring theory, but in the absence of additive inverses in semirings, their structure differs from that of ring theory. Due to this difference in both the theories, the role of $k-$ideals (if $x + y \in I, x \in I$, then $y \in I$) becomes significant in semirings. It is pertinent to note here that many results which are true for ideals in rings have been established, by many authors, for $k-$ideals in semirings (c.f. [5], [6], [7], [8]). This fact has motivated different researchers to settle the primary decomposition for k-ideals in semirings analogous to the primary decomposition theorem in rings (Lasker Noether theorem) which states that in a commutative Noetherian ring, every ideal can be described as a finite intersection of primary ideals.

The above result of ring theory is not true for arbitrary ideals in semirings as noticed in [1]. R. E. Atani and S. E. Atani first proved that in a commutative Noetherian semiring, every proper $k-$ideal can be represented as a finite intersection of $k-$primary ideals [1, Theorem 4]. As observed in [4], there are some errors in the results used to prove this theorem. For example, $I + Ra^n$ is not a $k-$ideal, even if $I$ is a $k-$ideal. But the authors of [1] took it for granted that the ideal $I + Ra^n$ is a $k-$ideal. P. Lescot [4] found these errors after observing in Example 6.2 that {0} ideal may not be a finite intersection of $k-$primary ideals in a commutative Noetherian semiring. With these observation, P. Lescot developed the theory of weak primary decomposition for semirings of characteristic 1. But still the question of settling the primary decomposition for a proper $k-$ideal (other than {0} and semiring $R$) remained unsolved. In this direction, S. Kar et. al. [3, Theorem 4.4] proved that every proper $k-$ ideal of a commutative Noetherian semiring $R$ can be expressed as a finite intersection of

$k$ −irreducible ideals, where an ideal $I$ is $k$ −irreducible if for any two $k$ −ideals $J, K$ of $R$, $I = J \cap K$ implies that either $I = J$ or $I = K$. They tried to prove the primary decomposition by observing that in a commutative Noetherian semiring, every $k$ −irreducible ideal is a $k$ −primary ideal [3, Theorem 4.6]. But the proof of this result has the same errors again and the problem still remained unsolved. In this paper, we provide a correct proof of Theorem 4.6 of [3], which settle the primary decomposition for $k$ −ideals in commutative Noetherian semirings. In order to establish the uniqueness of primary decomposition for semirings, first we prove some basic results required for uniqueness and finally also establish the uniqueness of primary decomposition for $k$ −ideals of commutative Noetherian semirings, that is, if $I = \cap_{i=1}^{n} Q_i$, where each $Q_i$ is a primary and $\sqrt{Q_i} = P_i$, then the set $\{P_1, P_2, \ldots P_n\}$ is independent of particular choice of decomposition of $I$.

## 2. Primary Decomposition of $k$ −Ideals in Semirings

In this section, we prove the primary decomposition and uniqueness theorem for $k$ −ideals in commutative Noetherian semirings. First, we recall some definitions and results from [2] and [3] which are necessary to prove the main results of this section.

**Definition 2.1 [2].** A semiring is a nonempty set $R$ on which operations of addition and multiplication have been defined such that for $r \in R$, the following conditions are satisfied:

(1) $(R, +)$ is a commutative monoid with identity element $0_R$;

(2) $(R, .)$ is a monoid with identity element $1_R$;

(3) Multiplication distributes over addition from either side;

(4) $0_R . r = 0_R = r . 0_R$;

(5) $1_R \neq 0_R$.

A semiring $R$ is said to be commutative, if $(R, .)$ is a commutative monoid. Throughout this paper, $R$ is a commutative semiring with identity.

**Definition 2.2 [2].** Let $I$ be an ideal of a semiring $R$. Then the radical of $I$, denoted by $\sqrt{I}$ is defined as $\sqrt{I} = \{a \in R | a^n \in I \text{ for some } n \geq 1\}$.

**Definition 2.3 [2].** A proper ideal $I$ of a semiring $R$ is called a prime ideal, if $ab \in I$ for $a, b \in R$ implies that either $a \in I$ or $b \in I$.

**Definition 2.4 [2].** A proper ideal $I$ of a semiring $R$ is said to be a primary ideal of $R$ if for any $x, y \in I, xy \in I$ implies that either $x \in I$ or $y \in \sqrt{I}$. If $\sqrt{I} = P$, then $I$ is called $P$ −primary.

**Definition 2.5 [3].** A proper $k$ −ideal $I$ of a semiring $R$ is called a $k$ −irreducible ideal of $R$ if for any two $k$ −ideals $J, K$ of $R$, $I = J \cap K$ implies that either $I = J$ or $I = K$.

We now present the notion of primary decomposition of k-ideals as follows:

**Definition 2.6.** Let $I$ be a $k$ −ideal of a semiring $R$. Then $I$ is said to have a primary decomposition if $I$ can be expressed as $I = \cap_{i=1}^{n} Q_i$, where each $Q_i$ is a primary ideal of $R$.

Also, a primary decomposition of the type $I = \bigcap_{i=1}^{n} Q_i$ with $\sqrt{Q_i} = P_i$, is called a reduced primary decomposition of $I$, if $P_i's$ are distinct and $I$ cannot be expressed as an intersection of a proper subset of ideals $Q_i$ in the primary decomposition of $I$. A reduced primary decomposition can be obtained from any primary decomposition by deleting those $Q_j$ that contains $\bigcap_{\substack{i=1 \\ i \neq j}}^{n} Q_i$ and grouping together all distinct $\sqrt{Q_i}'s$.

The main aim of this paper is to prove the existence and uniqueness of primary decomposition for $k-$ideals in a commutative Noetherian semiring. First, we prove the existence part as stated below:

**Theorem 2.7 (Primary Decomposition of $k-$Ideals).** Let $I$ be a $k-$ideal of a Noetherian semiring $R$. Then $I$ can be represented as a finite intersection of primary ideals of $R$, that is, $I$ has a reduced primary decomposition.

The decomposition of $k-$ideals in terms of $k-$irreducible ideals has already been proved in [3] as follows:

**Theorem 2.8 [3, Theorem 4.4].** Every proper $k-$ideal of a commutative Noetherian semiring $R$ can be expressed as a finite intersection of $k-$irreducible ideals.

Therefore, primary decomposition for $k-$ideals in semirings follows immediately if we prove

**Theorem 2.9.** Let $R$ be a commutative Noetherian semiring. Then every $k-$irreducible ideal of $R$ is a primary ideal of $R$.

Kar et. al. [3, Theorem 4.6] made an attempt to prove the same, but the proof of the result has some errors as mentioned in the introduction. Before we prove the required result, we analyse some common errors committed by many authors regarding $k-$ideals, and prove some facts about k-ideals. First, we note that any ideal $I = <a>$ is not a $k-$ideal in semirings as observed by J. S. Golan in [2, Example 6.17]. That is, if $I = <(1+x)>$ is an ideal of $\mathbb{N}[x]$ (semiring of polynomials over non-negative integers $\mathbb{N}$ in the indeterminate $x$), then $(1+x)^3 = (x^3 + 1) + 3x(1+x) \in I$, $3x(1+x) \in I$, but $(x^3+1) \notin I$ implies that $I$ is not a $k-$ideal of $\mathbb{N}[x]$. Any ideal generated by a single element becomes a $k-$ideal, if we impose some conditions on a semiring as shown below:

**Lemma 2.10.** Let $R$ be an additively cancellative, yoked and zerosumfree semiring. Then for any $a \in R$, the ideal $I = <a>$ is a $k-$ideal of $R$.

**Proof.** Let $a_1 + a_2$, $a_1 \in I$. Then there exist some $r_1$, $r_2 \in R$ such that $a_1 + a_2 = ar_1$ and $a_1 = ar_2$. As $R$ yoked, so there exists some $r_3 \in R$ such that $r_1 + r_3 = r_2$ or $r_2 + r_3 = r_1$. If $r_1 + r_3 = r_2$, then $ar_1 = ar_2 + a_2 = ar_1 + ar_3 + a_2$ implies that $ar_3 + a_2 = 0$, as $R$ is additively cancellative. Also, $a_2 = ar_3 = 0 \in I$, as $R$ is zerosumfree. Further, if $r_2 + r_3 = r_1$, then $ar_2 + a_2 = ar_1 = a(r_2 + r_3) = ar_2 + ar_3$ implies that $a_2 = ar_3 \in I$, since $R$ is additively cancellative. Thus, $I = <a>$ is a $k-$ideal of $R$. ∎

The semiring considered in above example is additively cancellative, zerosumfree, but it is not yoked, for let $f(x) = 5x^2 + 9x + 2$ and $g(x) = 11x^2 + 3x + 5$ be two polynomials in

$\mathbb{N}[x]$, then there exists no $h(x) \in \mathbb{N}[x]$ such that either $f(x) + h(x) = g(x)$ or $g(x) + h(x) = f(x)$.

Similar to an ideal generated by a single element, the sum of two $k$-ideals may not be a $k$-ideal in a semiring. There are plenty of $k$-ideals in the semiring $\mathbb{N}$, but their sum is not a $k$-ideal. However, the sum of two $k$-ideals is a $k$-ideal in a lattice ordered semiring (c.f. [2], Corollary 21.22). While proving Theorem 2.9, the authors wrongly used that the ideals $(Q+<a^k>)$ and $(Q+<b>)$ are $k$-ideals. In view of the above observations, Theorem 2.9 follows verbatim as proved in [3, Theorem 4.6] for additively cancellative, yoked and zerosumfree semirings, because in this case, both an ideal generated by a single element and sum of two $k$-ideals is a $k$-ideal. Here, we give a proof of Theorem 2.9 without resorting to these restrictions.

**Proof of Theorem 2.9.** Let $Q$ be a $k$-irreducible ideal of a Noetherian semiring $R$. Let $ab \in Q$ be such that $b \notin Q$. Now, we construct two ideals $I$ and $J$ of $R$ as follows : $I = <a^k> + Q$ and $J = <b> + Q$. Then, clearly, $Q \subseteq I \cap J$. Let $y \in I \cap J$. Then $y = a^n z + q$ for some $z \in R$ and $q \in Q$. Again $aJ \subseteq Q$ (since $ab \in Q$) and so $ay \in Q$ (since $y \in J$). Therefore, $ay = a^{n+1}z + aq$. Thus, we get that $a^{n+1}z + aq \in Q$. Also, $aq \in Q$, since $q \in Q$ and $Q$ is an ideal of $R$. Thus, it follows that $a^{n+1}z \in Q$, since $Q$ is a $k$-ideal of $R$. Construct a set $A_n = \{x \in R \mid a^n x \in Q\}$. It is easy to check that $A_n$ is an ideal of $R$ and $A_1 \subseteq A_2 \ldots$ is an ascending chain of ideals. Since $R$ is Noetherian, $A_n = A_{n+1} = \ldots$ for some $n \in \mathbb{Z}^+$. Again $a^{n+1}z \in Q$ implies that $z \in A_{n+1} = A_n$. It demonstrates that $a^n z \in Q$ which implies that $y \in Q$. Thus, $I \cap J = Q$.

Let $Jac(R)$ denote the Jacobson radical of semiring $R$, that is, the intersection of all maximal $k$-ideals of $R$. Assume that $A$ is an ideal of $R$ and $Jac(A)$ denotes the intersection of all maximal $k$-ideals of $R$ containing $A$.

If $\bar{A}$ denotes the $k$-closure (i.e. $\bar{A} = \{a \in A\} \mid a + b = c$ for some $b, c \in A$) of an ideal $A$ of $R$, we have $(\overline{\sqrt{I \cap J}}) = \sqrt{\bar{I}} \cap \sqrt{\bar{J}}$, where $\sqrt{P}$ denotes the nil radical of an ideal $P$ (c.f. [4], Lemma 2.2). We know that $\sqrt{P} \subseteq Jac(P)$. Thus, $\bar{I} \cap \bar{J} \subseteq \sqrt{\bar{I}} \cap \sqrt{\bar{J}} = \overline{\sqrt{I \cap J}} \subseteq Jac(\overline{I \cap J})$. Therefore, $Jac(\bar{I} \cap \bar{J}) \subseteq Jac(\overline{I \cap J})$. Again $\overline{I \cap J} \subseteq \bar{I} \cap \bar{J} \Rightarrow Jac(\overline{I \cap J}) \subseteq Jac(\bar{I} \cap \bar{J})$. Thus, $Jac(\overline{I \cap J}) = Jac(\bar{I} \cap \bar{J}) = Jac(\bar{I}) \cap Jac(\bar{J})$. So we have $Jac(\overline{I \cap J}) = Jac(\bar{I}) \cap Jac(\bar{J})$.

Now, $Q = I \cap J \Rightarrow \bar{Q} = \overline{I \cap J} \Rightarrow Q = \overline{I \cap J}$, since $Q$ is a $k$-ideal of $R$. This shows that $Jac(Q) = Jac(\overline{I \cap J}) = Jac(\bar{I}) \cap Jac(\bar{J})$, that is, $Jac(Q) = Jac(\bar{I}) \cap Jac(\bar{J})$.

Again $Q$ is a $k$-irreducible ideal, which implies that $Jac(Q)$ is $k$-irreducible. $Jac(Q) = Q \cap Jac(R)$ but $Jac(Q) \neq Jac(R)$. Thus, $Jac(Q) = Q$, since $Jac(Q)$ is $k$-irreducible. Accordingly, $Q = Jac(\bar{I}) \cap Jac(\bar{J})$, where each of $Jac(\bar{I})$ and $Jac(\bar{J})$ are $k$-ideals of $R$. Now $b \in J$ implies that $b \in Jac(\bar{J})$, but $b \notin Q$, that is, $Q \neq Jac(\bar{J})$. So $Q = Jac(\bar{I})$, since $Q$ is $k$-irreducible. Now $a^n \in I \Rightarrow a^n \in Jac(\bar{I}) = Q$. Hence $Q$ is $k$-primary. ∎

**Remark 2.11.** Now Theorem 2.7 follows by combining [3, Theorem 4.4] and Theorem 2.9. It is important to note here that all $Q_i's$ in the primary decomposition of a $k$-ideal $I$ have an additional property that these are also $k$-ideals. Now it only remains to prove the uniqueness of primary decomposition of a $k$-ideal $I$.

The following lemma will be used to prove the uniqueness of reduced primary decomposition of $k$-ideals in semirings.

**Lemma 2.12.** Let $R$ be a semiring and $Q$ a $P$ −primary ideal of $R$. Then we have

(i) If $x \in R$, $(Q:x) = \{r \in R \mid rx \in Q\}$ is an ideal of $R$;

(ii) If $x \in Q$, then $(Q:x) = R$;

(iii) If $x \notin P$, then $(Q:x) = Q$;

(iv) If $x \notin Q$, then $(Q:x)$ is $P$ −primary.

**Proof.** (i) Let $r_1, r_2 \in (Q:x)$ and $a \in R$. Then $r_1x, r_2x \in Q$ implies that $(r_1 + r_2)x \in Q$ and $r_1 ax \in Q$ as $Q$ is an ideal of $R$. Thus, $(Q:x)$ is an ideal of $R$.

(ii) If $x \in Q$, then $Rx \subseteq Q$, as $Q$ is an ideal of $R$ which implies that $R \subseteq (Q:x)$. Also, $(Q:x)$ is an ideal of R and so $(Q:x) = R$.

(iii) Assume that $x \notin P$. If $a \in Q$, then $ax \in Q$ implies that $a \in (Q:x)$. For the converse part, suppose that $b \notin Q$ and $xb \in Q$. Then, $x \in \sqrt{Q} = P$ as $Q$ is $P$ −primary ideal of $R$ which contradicts that $x \notin P$. Thus, $xb \notin Q$ and so $b \notin (Q:x)$.

(iv) Suppose that $x \notin Q$. If $y \in (Q:x)$, then $xy \in Q$ implies that $y \in \sqrt{Q} = P$ as $Q$ is $P$ −primary ideal of $R$. Thus, $Q \subseteq (Q:x) \subseteq P$ implies that $P = \sqrt{Q} \subseteq \sqrt{(Q:x)} \subseteq \sqrt{P}$. Also, $\sqrt{P} = P$, as $P$ is a prime ideal of $R$ and so $P = \sqrt{(Q:x)}$. Now, we show that $(Q:x)$ is a primary ideal of $R$. The ideal $(Q:x)$ is proper as $x \notin Q$ and so, $1 \notin (Q:x)$. Assume that $ab \in (Q:x)$ and $b \notin \sqrt{(Q:x)}$ for $a, b \in R$. Then, $abx \in Q$ and $Q$ is a $P$ −primary ideal of $R$ which implies that either $ax \in Q$ or $b \in P = \sqrt{(Q:x)}$. Thus, $a \in (Q:x)$ as $b \notin \sqrt{(Q:x)}$. ∎

We now prove the uniqueness of the reduced primary decomposition of a $k$ −ideal of a semiring as follows:

**Theorem 2.13 (Uniqueness of Primary Decomposition).** Let $R$ be a commutative Noetherian semiring and $I$ a $k$ −ideal of $R$. If $I = \cap_{i=1}^{n} Q_i$ is a reduced primary decomposition of $I$ with $\sqrt{Q_i} = P_i$ for $i = 1,2....n$, then $\{P_1, P_2.... P_n\} = \{$Prime ideals $P \mid \exists\ x \in R$ such that $P = \sqrt{(I:x)}\}$. The set $\{P_1, P_2.... P_n\}$ is independent of the particular reduced primary decomposition chosen for $I$.

**Proof.** Let $x \in R$. Then $\sqrt{(I:x)} = \sqrt{(\cap_{i=1}^{n} Q_i : x)} = \cap_{i=1}^{n} \sqrt{(Q_i:x)} = \cap_{\substack{i=1 \\ x \notin Q_i}}^{n} P_i$ by Lemma 2.12 (iv) and therefore, $\sqrt{(I:x)} \subseteq P_i$ for all $i = 1,2....n$. Also, if $\sqrt{(I:x)}$ is prime, then $\prod_{\substack{i=1 \\ x \notin Q_i}}^{n} P_i \subseteq \cap_{\substack{i=1 \\ x \notin Q_i}}^{n} P_i = \sqrt{(I:x)}$ implies that $P_i \subseteq \sqrt{(I:x)}$ for some $i = 1,2...n$. Thus, we have $\{$Prime ideals $P \mid \exists\ x \in R$ such that $P = \sqrt{(I:x)}\} \subseteq \{P_1, P_2.... P_n\}$.

On the other hand, for $i \in \{1,2....n\}$, we have $\cap_{\substack{j=1 \\ j \neq i}}^{n} Q_j \nsubseteq Q_i$, as the primary decomposition is reduced. So there exists some $x_i \in \cap_{\substack{j=1 \\ j \neq i}}^{n} Q_j$ and $x_i \in Q_i$. If $y \in (Q_i:x_i)$, then $yx_i \in Q_i$, and $yx_i \in (\cap_{\substack{j=1 \\ j \neq i}}^{n} Q_j) \cap Q_i = I$ which implies that $y \in (I:x_i)$,. Thus, $(Q_i:x_i) \subseteq (I:x_i) \subseteq (Q_i:x_i)$, as $I \subseteq Q_i$. So $(Q_i : x_i) = (I : x_i)$ implies that $\sqrt{(Q_i : x_i)} = \sqrt{(I : x_i)} = P_i$ by Lemma 2.12 (iv). Hence $\{P_1, P_2, \ldots, P_n\} = \{$Prime ideals $P \mid \exists\, x \in R$ such that $P = \sqrt{(I : x_i)}\}$. ∎